\documentclass{artjlt}

\usepackage{amsmath,amsfonts,amssymb}

\title{Relationship between  Nichols braided Lie algebras and Nichols  algebras} 
\author{Weicai Wu,  Shouchuan Zhang and  Yao-Zhong Zhang}                 
\lastname{Wu, Zhang and Zhang}  



\address{%
Weicai Wu\\               
Department  of Mathematics\\ 
Hunan University\\
Changsha  410082\\ 
P.R. China \\            
weicaiwu@hnu.edu.cn                
}
\address{%
Shouchuan Zhang\\               
Department  of Mathematics\\ 
Hunan University\\
Changsha  410082\\ 
P.R. China \\            
sczhang@hnu.edu.cn               
}
\address{%
Yao-Zhong Zhang\\               
School of Mathematics and Physics\\ 
The University of Queensland\\
Brisbane 4072\\ 
Australia \\            
yzz@maths.uq.edu.au               
}

\begin{document}


\maketitle

\begin{abstract}  We establish the relationship among Nichols algebras, Nichols braided Lie algebras
 and  Nichols Lie algebras. We prove two results: (i) Nichols algebra $\mathfrak B(V)$ is
finite-dimensional if and only if Nichols braided Lie algebra $\mathfrak L(V)$ is finite-dimensional
if there does not exist any $m$-infinity element in $\mathfrak B(V)$; (ii) Nichols Lie algebra
 $\mathfrak L^-(V)$ is infinite  dimensional if $ D^-$ is infinite. We give the sufficient  conditions 
for Nichols braided Lie algebra $\mathfrak L(V)$ to be a homomorphic  image of  a braided Lie algebra
generated by  $V$ with defining relations.
\end {abstract}

\numberwithin{equation}{section}

\section{Introduction}\label {s0}
The theory of Lie superalgebras has been developed systematically, which includes the
representation theory and classifications of simple Lie superalgebras and their varieties \cite {Ka77}
. In many physical applications or in pure mathematical interest, one has to consider
not only $\mathbb Z_2$-- or $\mathbb Z$-- grading but also G-grading of Lie algebras, where G is an abelian
group equipped with a skew symmetric bilinear form given by a 2-cocycle. Lie algebras in symmetric
and more general categories were discussed in \cite {GRR95, Gu86, ZZ04}. A sophisticated
multilinear version of the Lie bracket was considered in \cite {Kh99a, Pa98}. Various generalized
Lie algebras have already appeared under different names, e.g. Lie color algebras, $\epsilon $ Lie
algebras \cite {Sc79}, quantum and braided Lie algebras \cite {Ma94, KS97},
generalized Lie algebras \cite {BFM96} and H-Lie algebras \cite {BFM01}.
In \cite {Ar11}, a Milnor--Moore type theorem for primitively generated braided bialgebras
was obtained by means of braided Lie algebras.
The question of finite-dimensionality of Nichols algebras dominates an important part of the recent
 developments in the theory of (pointed) Hopf algebras.
 The interest in this problem comes from the Lifting Method by Andruskiewitsch and Schneider to classify
 finite (Gelfand-Kirillov) dimensional pointed Hopf algebras, which are generalizations of quantized
 enveloping algebras of semi-simple Lie algebras. The  classification of finite dimensional pointed 
Hopf algebras  was studied in \cite {AS02, AHS08,AS10, He05, He06a, He06b, WZZ}.

This paper provides a new method to determine whether a Nichols algebra is finite dimensional or not.


Let $\mathfrak B(V)$ be the Nichols algebra of vector space $V$.
Let  $\mathfrak L(V)$ , $\mathfrak L^-(V)$ and $\mathfrak L_c(V)$ denote the  braided Lie algebras
 generated by $V$ in $\mathfrak B(V)$ under Lie operations $[x,  y] = yx - p_{yx}xy$,
$[x,  y]^- = xy - yx$ and $[x,  y]_c = xy - p_{xy}yx$, respectively, for any homogeneous elements
$x,  y \in \mathfrak B(V)$.   $(\mathfrak L(V),  [\ ])$,  $(\mathfrak L^-(V),  [\ ]^-)$ and
$(\mathfrak L_c(V),  [\ ]_c)$ are called Nichols braided Lie algebra,  Nichols Lie algebra
and Nichols braided m-Lie algebra of $V$,  respectively.
It is clear that $(\mathfrak L(V),  [ \ ])$ and $(\mathfrak L_c(V),  [ \ ]_c)$
are equivalent as vector spaces. If $\mathfrak B(V)$ is finite dimensional
then $\mathfrak B(V)$ is nilpotent,  so
$(\mathfrak L(V),  [\ ])$ and $(\mathfrak L^-(V),  [\ ]^-)$ also are nilpotent.

In this paper we prove the following two results: (i) $\mathfrak B(V)$ is
finite-dimensional if and only if $\mathfrak L(V)$ is finite-dimensional when there does not
exist any $m$-infinity element; (ii) $\mathfrak L^-(V)$ is infinite dimensional
if $  D^-$ is infinite. We give the sufficient conditions for Nichols braided Lie algebra $\mathfrak L(V)$ 
to be a homomorphic image of a braided Lie algebra generated by  $V$ with defining relations.

This paper is organized as follows. In section \ref{s1} we recall some results on Nichols algebras
and fix the notation.
In section \ref{s2} we show that $\mathfrak L^-(V)$ is infinite
 dimensional if $ D^-$ is infinite.
In section \ref{s3} we prove that $\mathfrak B(V)$ is finite-dimensional if and only if
$\mathfrak L(V)$ is finite-dimensional when there does not exist any $m$-infinity element
in $\mathfrak B(V)$. In section \ref{s4} we present the condition for
$\mathfrak B(V) = F\oplus \mathfrak L(V).$ In section \ref{s5}
we give the sufficient conditions for Nichols braided Lie algebra $\mathfrak L(V)$ to be a homomorphic 
image of  a braided Lie algebra generated by  $V$ with defining relations.

Throughout, $\mathbb Z =: \{x | x \hbox { is an  integer }\}.$
$\mathbb N_0 =: \{x | x \in \mathbb Z,  x\ge 0\}.$
$\mathbb N =: \{x | x \in \mathbb Z,  x>0\}.$  $F$ denotes the base field of characteristic zero.

\section{ Preliminaries}\label {s1}
In this section we recall some results on Nichols algebras ( see \cite {AHS08}).
\begin {Lemma} \label {1} (see \cite {AHS08})If $(V,   \alpha,   \delta )$ is a $FG$-{\rm YD} module,
then tensor algebra $T(V)$ over $V$ is a $FG$-${\rm YD}$ module.
\end {Lemma}

\vskip.1in
If   $\{x_1,   \cdots,   x_n\}$ is  a basis of  vector space  $V$ and
$C(x_i\otimes  x_j) = q_{ij} x_j\otimes x_i$ with $q_{ij} \in F$,
then $V$  is called a braided vector space of diagonal type, $\{x_1,   \cdots,   x_n\}$
is called canonical basis and $(q_{ij})_{n\times n}$  is called braided matrix.
Throughout this paper all of braided vector spaces are connected and  of diagonal type
without special announcement. Let $G= \mathbb Z^n$ and
$E=\{e_{1}, e_{2}, \cdots, e_{n}\}, e_i =: ( \stackrel{i} { \overbrace{0,  \cdots 0,  1}},  \cdots 0)\in G$ ,  $1\le i \le n.$
Let $\chi$ be a bicharacter of $G$ such that  $\chi (e_i,  e_j) = p_{ij}$ and
$C(x_i\otimes  x_j) = \chi (e_i,  e_j) x_j\otimes x_j$. Let $S_{m}\in End_{k}(T(V)^{m})$
and $S_{1, j}\in End_{k}(T(V)^{j+1})$ denote the maps $S_{m}=\prod \limits_{j=1}^{m-1}(id^{\bigotimes m-j-1}\bigotimes S_{1, j})$ ,
 $S_{1, j}=id+C_{12}^{-1}+C_{12}^{-1}C_{23}^{-1}+\cdots+C_{12}^{-1}C_{23}^{-1}\cdots C_{j, j+1}^{-1}$
(in leg notation) for $m\geq 2$ and $j\in \mathbb N$. Then the subspace
$S=\bigoplus \limits _{m=2}^{\infty} ker S_{m}$ of the tensor $T(V)=\bigoplus \limits _{m=0}^{\infty} T(V)^{\bigotimes m}$
is a two-sided ideal,  and algebra $\mathfrak B(V)=T(V)/S$ is termed the Nichols algebra associated to $(V, C)$.
Define  linear map $p$ from $\mathfrak B(V) \otimes \mathfrak B(V) $ to $F$ such that
$p(u\otimes v) = \chi (deg (u),  deg (v)), $ for any homogeneous element $u,  v \in \mathfrak B(V).$
For convenience,  $p(u\otimes v)$ is denoted by $p_{uv}$.
 Let $A =: \{x_1,   x_2,   \cdots,   x_n\}$ be an alphabet, $A^*$ the set of all of  words in $A$
and $A^+ =: A^*\setminus 1$.  Define $x_1 < x_2< \cdots < x_n$ and the order on $A^*$ is the lexicographic ordering.
For the concept of words refer \cite {Lo83}. Let $| u| $ denote the length of word $u$.

\begin {Definition} \label {2} (\cite [Def. 1] {Kh99b})
A word $u$ is called a Lyndon word if $| u| =1$ or $| u| \geq2$,  and for each
representation $u=u_1u_2$,  where $u_1$and $u_2$ are nonempty
words,  the inequality $u<u_2u_1$ holds.
\end {Definition}
 Any word $u\in A^*$ has a unique decomposition into the product of
non-increasing sequence of Lyndon words by \cite [Th.5.1.5] {Lo83}.
If $u$ is a Lyndon word with $| u | >1$,  then there uniquely exist two Lyndon words $v$ and $w$
such that $u =vw$ and $v$ is shortest (see \cite [Prop. 5.1.3]{Lo83} and \cite {He07})
(the composition is called the Shirshov decomposition of $u$).

\begin {Definition} \label {3} We inductively define a linear map $[\  \ ]$ from $A^+$
to $\mathfrak B(V)$  as follows:
(1) $[u] =: u$  when  $u$ is a letter;
(2) $[u] =: [w][v]-p_{wv}[v][w]$ when $u$ is a Lyndon word with $| u | >1$ and $u =vw$ is a
Shirshov decomposition; (3) $[u] =: [[[l_1,  l_2],  l_3] \cdots l_m] $,
when $u = l_1l_2 \cdots l_m$ is a non-increasing product of Lyndon words, 
i.e. $l_1 \ge  l_2 \ge l_3 \ge \cdots \ge l_m $, and $l_i$   is a Lyndon word  for any $1\le i \le  m$.
\end {Definition}
Similarly, we inductively define a linear map $[\  \ ]^-$ from $A^+$ to $\mathfrak B(V)$  as follows:
(1) $[u]^-  = : u$  when  $u$ is a letter;
(2) $[u]^-  = : [w]^-[v]^- -[v]^-[w]^-$ when $u$ is a Lyndon word with $| u | >1$ and $u  = vw$ is a Shirshov decomposition;
(3) $[u]^-  = : [[[l_1,  l_2]^-,  l_3]^- \cdots l_m]^- $,  when $u  =  l_1l_2 \cdots l_m$ is a non-increasing product of Lyndon words.

$[u]$ is called a nonassociative word for any $u \in A^+$,  $[u]$ is called a standard
nonassociative word if $u$ is a Lyndon word. Every standard nonassociative word  is also
called  a super-letter.

\begin {Definition} \label {5} (\cite [Def. 6] {Kh99b})
A super-letter $[u]$ is said to be hard if it is not a linear
combination of products $[u_1][u_2]\cdot\cdot\cdot[u_i],  i\in
\mathbb N$,  where $[u_j]$ are super-letter with $[u]<[u_j]$,  $1\le j \le i$.
\end {Definition}

\begin {Definition} \label {6} (\cite [Def. 7]{Kh99b} or \cite [before Th. 10] {He07})
We say that the height of a super-letter $[u]$ with degree $d$ equals a natural number $h$
if $h$ is least with the following properties:

(1) $p_{uu}$ is a primitive root of unity of degree $t\ge 1, $ and $h=t$ ;

(2) super-word $[u]^h$ is a linear combination of super-words of degree $hd$ in greater
super-letters than $[u].$

If the number $h$ with above properties does not exist then we say that the height of $[u]$ is infinite.
\end {Definition}

Let $h_u$ denote the height of $u$.
Let ${\rm ord } (p_{uu})$ denote the  order of $p_{uu}$ with respect to multiplication.
$D =: \{[u] \mid [u] \hbox { is a hard super-letter }\}$. If $[u] \in D$ and
${\rm ord }  (p_{u, u}) =m>1$ with $h_u = \infty$, then $[u]$ is called an $m$-infinity element.
$P = :\{[u_{1}]^{k_1}[u_{2}]^{k_2}\cdots  [u_{s}]^{k_s}\  \mid \  [u_{i}]\in D,
k_i, s \in \mathbb N_0; 0 \le k_i < h_{u_i};  1\le i \le s;   u_s<u_{s-1}<\cdots< u_1\}$.
$\Delta ^+(\mathfrak B(V)): =  \{ \deg (u) \mid [u]\in D\}$.
$ \Delta (\mathfrak B(V)) := \Delta ^+(\mathfrak B(V)) \cup \Delta ^-(\mathfrak B(V))$,
which is called the root system of $V.$ If $ \Delta (\mathfrak B(V))$ is finite,
then it is called an arithmetic root system. Let
$E_e' =: \frac {1}{2} \mid \{ ([u], [v])\in D\times D  \mid [u], [v] \in D, p_{u, v} p_{v, u} \not=1\}\mid $.
Let $D^- : = \{ [u]^- \mid [u] \in D\}$ and ${\deg } (D^-) : = \{ \deg ([u]^- ) \mid [u] \in D\}$.
Let $E_e$ denote  the number of edges of generalized Dynkin diagram.
If $u=vw$ is the shirshov decomposition of $[u] \in D$, then $[v]$, $[w] \in D$, which are called  sons of $u$. If
$[u_1], [u_2], \cdots, [u_m]\in D$ and $u_{i+1}$  is  a son of $u_i$ for $1\le i \le m-1$,
then $u_2, u_3, \cdots, u_m$ are called descendants of $u_1.$

\begin{Remark} There does not exist  any $m$-infinity element in $\mathfrak B(V)$
if and only if   Property (P) in \cite [Section 2.2] {He05} holds.
\end{Remark}

\begin {Theorem} \label {8} (\cite [Th. 2] {Kh99b} or \cite [Th. 10] {He07})
$P$ is a basis of $\mathfrak B (V)$.
\end {Theorem}

\section{Relationship between Nichols algebras and Nichols Lie algebras}\label {s2}
In this section it is proved that Nichols Lie  algebra $\mathfrak L^-(V)$ is infinite dimensional
if $ D^-$ is infinite.

\begin {Lemma} \label {63} Assume that  $u = vw$ is a Shirshov decomposition of $u$.
       If $[u]\in D$,    then $[[v],   [w]]^-\neq0$. Furthermore, if $[v], [w] \in \mathfrak L ^-(V)$ (e.g. $\mid u \mid = 2$), then
$[u]^- \not =  0.$
\end {Lemma}
\begin{Proof} If $[[v],   [w]]^- = 0$,    then $[v][w] = [w][v]$,
we know $[[v],   [w]] = [w][v]-p_{wv}[v][w] = (1-p_{wv})[w][v]$,    it contradicts
to $[u] \in P$ and $    [w][v]\in P$. 
\end{Proof}

\begin {Theorem} \label {63'}  (i)  $\dim \mathfrak L^-(V) \ge \mid \deg (D^-)\mid -1 \ge n + E_e -1$,
where $E_e$ is the number of edges in generalized Dynkin diagram of $V$.
(ii) If $D^-$ is infinite,  then  $ \dim  \mathfrak L^-(V) = \infty$.
\end {Theorem}
\begin{Proof}  If $u_1, u_2, \cdots, u_m$ in $D^- \setminus 0$ with different degrees, then  $u_1, u_2, \cdots, u_m$  is  linearly independent.
It is clear that $D^- \subseteq L^-(V).$  Consequently, $\dim L^-(V) \ge \mid \deg ( D^-)\mid -1 $. Obviously,
there exists a line between $x_i$ and $x_j$  if and only if  $[x_ix_j]\in D$  with $i<j$, which implies $\mid \deg (D^-)\mid -1 \ge n +   E_e  -1$. \end{Proof}

\section{Relationship between Nichols algebras and Nichols braided Lie algebras}\label {s3}
In this section it is proved that  $\mathfrak B(V)$ is
finite-dimensional if and only if $\mathfrak L(V)$ is finite-dimensional when there does not
exist any $m$-infinity element in $\mathfrak B(V)$. Let $l_u(v) := [u, v]$ and $r_u(v) := [v, u]$ for any $u, v\in \mathfrak B(V).$

\begin {Lemma} \label {7} If $[u]$ is a nonassociative word,   then $[u] \in \mathfrak L(V)$.
\end {Lemma}
\begin{Proof}  By the definition of nonassociative words,   we have  $[u] \in \mathfrak L(V)$. \end{Proof}

\begin{Remark} If $| D | =\infty$, then $\dim \mathfrak L(V)=\infty$. \end{Remark}
\begin {Lemma} \label {11} If $[u]$ is a nonassociative word with
$ t \in \mathbb N$,  $t \le ord (p_{uu})$,   then $[u]^t \in \mathfrak L(V)$.
\end {Lemma}

\begin{Proof} Let $l_{[u]}^{0}[u]=:[u],   l_{[u]}^{i}[u]=:[ [u], l_{[u]}^{i-1}[u]],   i\geq1$.
   Obviously $l_{[u]}^{i}[u] \in \mathfrak L(V)$.  It is clear $l_{[u]}^{1}[u]=[[u],
 [u]]=[u]^2-p_{uu}[u]^2=(1-p_{uu})[u]^2$. By means of  induction,   we obtain
$l_{[u]}^k[u]=(1-p_{uu})(1-p_{uu}^2)\cdot\cdot\cdot(1-p_{uu}^{k})[u]^{k+1},  \forall \  1<k\in
\mathbb N $.
 We have $(1-p_{uu})(1-p_{uu}^2)\cdot\cdot\cdot(1-p_{uu}^{t-1}) \not=0$
since  $t \le  ord (p_{uu})$,   which  implies $[u]^t \in \mathfrak L(V).$ \end{Proof}

\begin {Theorem} \label {9}
If there does not  exist any $m$-infinity element in $\mathfrak B(V)$  and $1<ord (p_{uu}) <\infty$ for any  $u \in D$, then the following conditions are equivalent:
(i) $\mathfrak B(V)$ is finite-dimensional;
(ii) $\mathfrak L(V)$ is finite-dimensional;
(iii) $\Delta (\mathfrak B(V))$ is an arithmetic root system.
\end {Theorem}

\begin{Proof} It follows from   \cite [Section 2.2]{{He05}}  that (i) and (iii) are equivalent.
$(i)\Longrightarrow (ii)$. Assume that $\mathfrak B(V)$ is
finite-dimensional. Since $\mathfrak L(V)\subseteq \mathfrak B(V)$,
we have that $\mathfrak L(V)$ are finite-dimensional.
$(ii)\Longrightarrow (i)$.  Assume that $\mathfrak L(V)$ is finite-dimensional.
By Lemma \ref {7},  $D \subseteq \mathfrak L(V).$  Obviously,  $D \subseteq P$.
Therefore $D$ is linearly independent and $| D | \le \dim \mathfrak L(V) < \infty,$
$h_{u} <\infty$ since $1<ord (p_{uu}) <\infty$ for  $[u] \in D$. It follows from  Theorem \ref {8}
that $\dim \mathfrak B(V) < \infty$. \end{Proof}

\begin {Proposition} \label {63'''} Assume that  $V$ is a  Cartan type  with generalized
Cartan matrix $(a_{ij})_{n\times n}$  and
$1< ord (p_{uu}) < \infty $ for any $[u] \in D$.  If there does not exist any $m$-infinity element in $\mathfrak B(V)$, then the following conditions are equivalent.
(i) $\mathfrak L(V)$ is finite dimensional;
(ii) $(a_{ij})_{n\times n}$ is a Cartan matrix;
(iii) $\dim \mathfrak B(V)< \infty.$
\end {Proposition}

\begin{Proof} It follows from  \cite [Th. 2.10.2]{He05}, Theorem \ref {9} and Lemma  \ref {7}. \end{Proof}

\begin {Proposition} \label {16} If there exists  $ [u]\in D$ such that ${\rm ord }  (p_{uu})=\infty$,
then $\dim \mathfrak B(V)=\infty$ and $\dim \mathfrak L(V)=\infty$.
\end {Proposition}

\begin{Proof}  By Theorem \ref {8}, $\dim  \mathfrak B(V)=\infty$. By Lemma
\ref {11}, $\dim \mathfrak L(V)=\infty$.  \end{Proof}

\begin {Proposition} \label {21} If there exist $[u],[v]$
such that $p_{uu}^{i}p_{uv}p_{vu}\neq 1$ for $\ \forall \ 0\leq i\leq 2k-2$, $\ \forall \ k\in \mathbb N $,
then $[v][u]^k,[u][v][u]^{k-1},\ldots,[u]^{k}[v]\in \mathfrak L(V)$.
\end {Proposition}

\begin{Proof}  We first show
\begin {eqnarray}\label {e5}   &&p_{uu}^{k-t+i}p_{uv}r_{[u]}^{t-i-1}([u]^{i}l_{[u]}^{k-t+1}[v])+r_{[u]}^{t-i}([u]^{i}l_{[u]}^{k-t}[v])
    \nonumber \\
&&~~~~~~~~~~~~~~~~~~~~=   (1-p_{uu}^{2(k-t)+i}p_{uv}p_{vu})r_{[u]}^{t-i-1}([u]^{i+1}l_{[u]}^{k-t}[v])
\end {eqnarray}
for $\ \forall \ 1\leq t \leq k$ , $\ \forall \ 0\leq i \leq t-1$. In fact,
\begin{eqnarray*}
{\rm left~ hand~ side ~ of~ } (\ref {e5}) &=& p_{uu}^{k-t+i}p_{uv}r_{[u]}^{t-i-1}([u]^{i}[[u],
l_{[u]}^{k-t}[v]])+r_{[u]}^{t-i-1}([[u]^{i}l_{[u]}^{k-t}[v], [u]])\\
&=&r_{[u]}^{t-i-1}\left( p_{uu}^{k-t+i}p_{uv}[u]^{i}l_{[u]}^{k-t}[v][u]
          -p_{uu}^{k-t+i}p_{uv}[u]^{i}l_{[u]}^{k-t}[v][u]\right. \\
& &\left.-p_{uu}^{k-t+i}p_{uv}p_{uu}^{k-t}p_{vu}[u]^{i+1}l_{[u]}^{k-t}[v]+[u]^{i+1}l_{[u]}^{k-t}[v]
   \right)\\
&=&   (1-p_{uu}^{2(k-t)+i}p_{uv}p_{vu})r_{[u]}^{t-i-1}([u]^{i+1}l_{[u]}^{k-t}[v])\\
&=& {\rm right~ hand~ side~  of ~} (\ref {e5}).
\end {eqnarray*}

Let $B^{(i)}=(b^{(i)}_{rs})_{(k-i+1)\times (k-i+1)}$ be real matrices such that

  $ \left(\begin{array}{ccccccc}
              [u]^{i}l_{[u]}^{k-i}[v] \\
              r_{[u]}^{1}([u]^{i}l_{[u]}^{k-i-1}[v])\\
              r_{[u]}^{2}([u]^{i}l_{[u]}^{k-i-2}[v]) \\
              \cdots\\
              r_{[u]}^{k-1}([u]^{i}[v]) \\
                   \end{array}\right) $  {}   $ $=
 $ B^{(i)} ${}
 $\left( \begin{array}{c}
                [u]^{i}[v][u]^{k-i}\\

                [u]^{i+1}[v][u]^{k-i-1}\\

                [u]^{i+2}[v][u]^{k-i-2}\\
                  \cdots \\

                  [u]^{k}[v] \\
                  \end{array}\right)$ {} ${}
            $  and  $b_{11}^{(i)}=1$ for $0\le i \le t-1$.

We know
$| B^{(i)}| =\prod \limits _{t=i+1}^{k}(1-p_{uu}^{2(k-t)+i}p_{uv}p_{vu})| B^{(i+1)}| $ by (\ref {e5}).
Consequently,  $| B^{(0)}| =\prod \limits _{i=0}^{k-1}\prod \limits _{t=i+1}^{k}(1-p_{uu}^{2(k-t)+i}p_{uv}p_{vu})$
and $| B^{(0)}| =0$ is equivalent to $\prod \limits _{i=0}^{2k-2}(1-p_{uu}^{i}p_{uv}p_{vu})=0$.
This completes the proof. \end{Proof}

\vskip.1in
According to the above Proposition, we obtain immediately,
\begin {Corollary} \label {22} If there exist $[u], [v] \in D$
such that $p_{u,  u}=1$ and $p_{uv}p_{vu}\neq 1$,
then $\dim \mathfrak B(V)=\infty, \  \dim \mathfrak L(V)= \infty$ .
\end {Corollary}

\begin {Corollary} \label {23} If  $[u], [v]\in D,$  such that $p_{u,  u}=1$ and
 $p_{uv}p_{vu}=1$,  then $d_{1}d_{2}\cdots d_{k}[v]\in Fl_{[u]}^{k}[v]$ for $ \forall \ d_{i}=l_{[u]}$
or $r_{[u]}, 1\leq i \leq k$.
\end {Corollary}

\begin{Proof} We know
\begin {eqnarray}\label {e7}
p_{uv}r_{[u]}^{s}l_{[u]}^{k-s}[v]+r_{[u]}^{s+1}l_{[u]}^{k-s-1}[v]=(1-p_{uv}p_{vu})[u]r_{[u]}^{s}l_{[u]}^{k-s-1}[v]
\end {eqnarray}
for $\ \forall \ 0\leq s < k$ by Definition \ref {5}. Then
$l_{[u]}^{k}[v]=-p_{vu}r_{[u]}^{1}l_{[u]}^{k-1}[v]=\ldots=(-p_{vu})^{k}r_{[u]}^{k}[v]. $
On the other hand,  $r_{[u]}l_{[u]}[v]=l_{[u]}r_{[u]}[v]$.  This proves the corollary. \end{Proof}

\begin {Proposition} \label {41} If there exists $[u]\in D$
such that $p_{uv}p_{vu}=1$ for $  \ \forall  v \in D $ with $ [u]\neq [v]$,  then $\Delta(\mathfrak B(V)) $ is not an arithmetic root system
while $\mathfrak B(V) $ is connected Nichols algebra of diagonal type with $\dim V >1$. Moreover,
$\dim \mathfrak B(V)=\infty, \  \dim \mathfrak L(V)= \infty$.
\end {Proposition}

\begin{Proof} There exists a basis $\pi$ of $\Delta ( \mathfrak B (V))$ such that $\deg u \in \pi.$ Since  $\mathfrak B(V) $ is connected Nichols algebra and $n = \mid \pi \mid >1,$ there exists $\beta \in \pi\setminus \{\deg u\}$ such that $\chi (\deg u, \beta) \chi (\beta, \deg u) \not=1.$ By the definition of $\Delta ( \mathfrak B (V))$
there exists $[v] \in D\setminus \{u\}$ with $\deg v \in \{ \beta, - \beta\}.$  This yields a contradiction to $p_{uv}p_{vu} =1.$
\end{Proof}

\begin {Theorem} \label {51} If $\mathfrak B(V) $ is connected Nichols algebra of diagonal type with $\dim V>1$
and there does not exist any $m$-infinity elements,  then $\mathfrak B(V)$ is
finite-dimensional if and only if $\mathfrak L(V)$ is finite-dimensional.
\end {Theorem}

\begin{Proof} It follows from  Proposition \ref {9}, Proposition \ref {16} ,
 Corollary \ref {22} and Proposition \ref {41}. \end{Proof}

\vskip.1in
By \cite {ZZ04}, $(\mathfrak B(V),     [\ ]_c) $ is a braided  m-Lie algebra and
we have the braided Jacobi identity as follows:
\begin {eqnarray}\label {e2}
[ [u,   v],   w]= [u,    [v,   w]]  +p_{vw}^{ -1} [ [u,   w],   v]
+(p_{wv} -p_{vw}^{ - 1}) v\cdot [u,   w].
\end {eqnarray}

\begin {Lemma}  \label {12}
If $u$ and $ v$ are homogeneous elements in $\mathfrak L(V)$ with  $p_{uv}p_{vu}\not =  1$,
then $uv,   vu \in \mathfrak L(V)$. Furthermore, if  $u, v \in \mathfrak L (V)$, then $ [u, v]^- \in \mathfrak L(V).$
\end {Lemma}

\begin{Proof} $[u,   v] =  vu -p_{v,   u} uv$ and   $[v,   u] =  uv -p_{u,   v} vu$,   which implies that
$uv$ and $vu$ are a linear  combination of $[u,   v]$ and $[v,   u]$. \end{Proof}

\begin {Proposition} \label {3.33}  $\dim \mathfrak L(V) \ge \sum _{[u] \in D} (h_u -1) + E_e'.  $
\end {Proposition}
\begin{Proof} It follows from  Lemma \ref {11} and Lemma \ref {7}. \end{Proof}

\vskip.1in
Recall the dual $\mathfrak B(V^*) $ of Nichols algebra $\mathfrak B(V) $ of rank $n$
in \cite [Section 1.3]{He05} and \cite{He06b}. Let $y_{i}$ be a dual basis of
$x_{i}$. $\delta (y_i)  = g_i ^{-1} \otimes y_i$,    $g_i \cdot y_j = p_{ij}^{-1} y_j $
and $\Delta (y_i)  = g_i ^{-1} \otimes y_i +y_i \otimes 1.$ There exists a bilinear
 map $<,  >$ from $(\mathfrak B(V^*)  \# kG)  \times \mathfrak B(V) $ to $\mathfrak B(V) $
such that     $<y_i,    uv> = <y_i,    u>v +g_i^{-1}.u<y_i,    v>$ and
$<y_i,    <y_j,    u>> = <y_iy_j,    u>$  for any $u,    v\in  \mathfrak B(V)  $.
Furthermore,    for any $u\in \oplus _{i=1}^\infty \mathfrak B(V)_{(i)}$,
one has that  $u=0$ if and only if $<y_i,    u> = 0$ for any $1\leq i \leq n.$
 Let $i$ denote $x_{i}$ in short, sometimes.

\begin{Lemma}\label{14} Let $l_{i}^{0}[j]=[j]$ ,  $l_{i}^{k}[j]=[i,  l_{i}^{k-1}[j]]$ ,
$r_{i}^{0}[j]=[j]$ , $r_{i}^{k}[j]=[r_{i}^{k-1}[j], i]$ ,  $k\geq 1$ . Then we have

(i)   $<y_j,    l_{i}^{k}[  j]>=0,  <y_i,    r_{i}^{k}[  j]>=0,  \forall\ k\geq 1$;

(ii) the following conditions are equivalent:
 (1) $l_{i}^{k}[j]=0$ ;
    (2) $r_{i}^{k}[j]=0$ ;
    (3)  $(k)_{p_{ii}}^{!}\prod \limits _{t=0}^{k-1}(p_{ii}^{t}p_{ji}p_{ij}-1)=0$  .
\end {Lemma}

\begin{Proof}
(i) It is clear $<y_j,    l_{i}^{k}[  j]>=0$,  $ <y_i,    r_{i}^{k}[  j]>=<y_i,
ir_{i}^{k-1}[  j]-p_{ii}^{k-1}p_{ij}r_{i}^{k-1}[  j]i>
=r_{i}^{k-1}[  j]-p_{ii}^{k-1}p_{ij}p_{ii}^{-(k-1)}p_{ij}^{-1}r_{i}^{k-1}[  j]=0$.

(ii) By means of induction, we obtain
$<y_{i}, l_{i}^{k}[j]>=p_{ii}^{-(k-1)}p_{ij}^{-1}(1-p_{ii}^{k-1}p_{ij}p_{ji})(1+p_{ii}
+\cdots p_{ii}^{k-1})l_{i}^{k-1}[j]$,  then
$<y_{j}y_{i}^k, l_{i}^{k}[j]>=p_{ii}^{-\sum \limits_{i=1}^{k-1}i}p_{ij}^{-k}\prod
\limits _{t=0}^{k-1}(1-p_{ii}^{t}p_{ij}p_{ji})(1+p_{ii}+\cdots p_{ii}^{t})$,
(1) is equivalent to (3) by (i). On the other hand,
$<y_{j}, r_{i}^{k}[j]>=p_{ji}^{-k}\prod \limits _{t=0}^{k-1}(1-p_{ii}^{t}p_{ij}p_{ji})[i]^k$,
$<y_{i}^{k}y_{j}, r_{i}^{k}[j]>=p_{ii}^{-\sum \limits_{i=1}^{k-1}i}p_{ji}^{-k}\prod
\limits _{t=0}^{k-1}(1-p_{ii}^{t}p_{ij}p_{ji})(1+p_{ii}+\cdots p_{ii}^{t})$.
One knows that (2) is equivalent to (3) by (i). This proves the lemma.
\end{Proof}

\section{Conditions for $\mathfrak B(V) = F\oplus \mathfrak L(V).$}\label {s4}
In this section we give  the sufficient conditions for  $\mathfrak B(V) = F\oplus \mathfrak L(V).$

\begin{Lemma}\label{13} (\cite [Lemma 3.1] {WZZ})
 {(\rm i)} If $| u | = | v |$,  then $u<v$ if and only if $uw < vw.$
 {(\rm ii)} If $u=vw$ is the Shirshov decomposition of Lyndon word $u$  and  $[u]$ is  hard,
then both  $[v]$ and $[w]$ are  hard too.
 \end{Lemma}

\begin {Lemma} \label {10} If there exist $x_{i}, x_{j}, i\neq j$ such that $p_{ij}p_{ji}=1$,
then $\mathfrak B(V) \neq F\oplus \mathfrak L(V).$
\end {Lemma}

\begin{Proof}  It is clear $[x_{i}, x_{j}]=[x_{j}, x_{i}]=0$ and $x_{i}x_{j}=p_{ij}x_{j}x_{i}\neq 0$.
Then $x_{i}x_{j}$ or $x_{j}x_{i}\in P$ and $x_{i}x_{j}, x_{j}x_{i}\notin \mathfrak L(V).$  \end{Proof}

\begin {Corollary} \label {60}
If $\mathfrak B(V) $ is connected Nichols algebra of rank $>3$ of diagonal type and $\Delta(\mathfrak B(V)) $
is arithmetic root systems, then $\mathfrak B(V) \neq F\oplus \mathfrak L(V).$
\end {Corollary}

\begin{Proof} It is clear from \cite [Table A.1] {He05},  \cite [Table A.2] {He05},
\cite [Table B] {He06a} and \cite [Table C] {He06a}. \end{Proof}

\begin {Example}  \label {50}
If  $\begin{picture}(100,       20)
\put(10,      1) {\makebox(0, 0) [t]{$\bullet$}}
\put(80,      1){\makebox(0, 0) [t]{$\bullet$}}
 \put(10,     -1) {\line(1,  0) {70}}
 \put(5,       9) {$\zeta$}
 \put(40,       4) {$-\zeta$}
 \put(75,       9) {$-1$}
  \put(90,       1) {$,    $}\ \ \ \put(95,       1)  {$ \zeta \in R_3$,}
\end{picture}$ \\
then  $D=\{   [x_1],      [x_2],    [x_1,    x_2],      [x_1,  [x_1,    x_2]]\}, \dim \mathfrak B(V) =2^{2}3^{2}$  and
$\mathfrak B(V) = F\oplus \mathfrak L(V)$.
\end {Example}

\begin{Proof}
  Assume that $  [u]$ is a hard super-letter or zero and  $u=vw$ is the Shirshov decomposition of $u$ when $[u] \not=0$.
 We show  $  [u] \in D$ step by step for the length $| u | $ of $u$.

{\rm (a)}   $| u| =2$,     then $u=  [1,   2]$ by Lemma \ref {14}.

{\rm (b)}  $| u| =3$,    then $  [u]=
  [1,   [1,   2]] $ and $[[1 , 2],  2]=0$ by Lemma \ref {14}.

{\rm (c)}  $| u| =4$. then $
  [1,  [1,   [1,   2]] =0$ by Lemma \ref {14}.

{\rm (d)} $| u| =5$, then
\begin{eqnarray*}
<y_1,    [[1,   [1,   2]],     [1,   2]]>&=&<y_1,   [1,   2][1,   [1,   2]]-p_{11}^{2}p_{12}p_{21}^{2}p_{22}[1,   [1,   2]][1,   2]>\\
&=&p_{12}^{ -1}(1-p_{12}p_{21})2[1,   [1,   2]]\\
& & +p_{11}^{-1}p_{12}^{-1}[1,   2]p_{11}^{-1}p_{12}^{-1}(1-p_{11}p_{12}p_{21})(1+p_{11})[1,   2]\\
& & - p_{11}^{2}p_{12}p_{21}^{2}p_{22}p_{11}^{-1}p_{12}^{-1}(1-p_{11}p_{12}p_{21})
(1+p_{11})[1,   2][1,   2] \\
& & -p_{11}^{2}p_{12}p_{21}^{2}p_{22} p_{11}^{-2}p_{12}^{-1}[1,[1, 2]]p_{12}^{ -1}(1-p_{12}p_{21})2\\
&=& p_{12}^{ -1}(1-p_{12}p_{21})[[1,   [1,   2]], 2]\\
& & +(p_{11}^{-1}p_{12}^{-1}-p_{11}^{2}p_{12}p_{21}^{2}p_{22})p_{11}^{-1}p_{12}^{-1}\\
& &~~~\times (1-p_{11}p_{12}p_{21})(1+p_{11})[1,   2]^2\\
&=&p_{12}^{ -1}(1-p_{12}p_{21})\left\{p_{12}^{-1}p_{22}^{-1}[[1, 2],[1,  2]]\right.\\
& &~~ \left.+(p_{21}p_{22}-p_{12}^{-1}p_{22}^{-1})[1,   2]^2\right\}\\
& &+(1-p_{11}^{3}p_{12}^{2}p_{21}^{2}p_{22})p_{11}^{-2}p_{12}^{-2}(1-p_{11}p_{12}p_{21})
(1+p_{11})[1,   2]^2\\
&=&p_{12}^{ -1}(1-p_{12}p_{21})(p_{21}p_{22}-p_{11}p_{21})[1,   2]^2\\
& &+(1-p_{11}^{3}p_{12}^{2}p_{21}^{2}p_{22})p_{12}^{-2}\\
& &~~~\times (p_{11}^{-2}+p_{11}^{-1}-p_{11}^{-1}p_{12}p_{21}-p_{12}p_{21})[1,   2]^2\\
&=&p_{12}^{ -2}\{p_{12}p_{21}(1-p_{12}p_{21})(p_{22}-p_{11})\\
& &+(1-p_{11}^{3}p_{12}^{2}p_{21}^{2}p_{22})\\
& &~~~\times (p_{11}^{-2}+p_{11}^{-1}-p_{11}^{-1}p_{12}p_{21}-p_{12}p_{21})\}[1,   2]^2\\
&=&p_{12}^{ -2}\left(p_{11}^{ -2}(1-p_{11}p_{12}p_{21})(1-p_{11}^{2}
p_{12}p_{21})(1+p_{11}^{ 2}p_{12}p_{21}p_{22})\right.\\
& &\left.+p_{11}^{-1}(1-p_{11}^{ 2}p_{12}p_{21})(1-p_{11}p_{12}^{ 2}p_{21}^{ 2}p_{22})\right)[1,   2]^2\\
&=&p_{12}^{ -2}(1-p_{11}^{2}p_{12}p_{21})\{p_{11}^{ -2}(1-p_{11}p_{12}p_{21})
(1+p_{11}^{ 2}p_{12}p_{21}p_{22})\\
& & +p_{11}^{-1}(1-p_{11}p_{12}^{ 2}p_{21}^{ 2}p_{22})\}[1,   2]^2\\
&=&p_{12}^{ -2}(1-p_{11}^{2}p_{12}p_{21})\{p_{11}^{ -2}(1-p_{11}p_{12}p_{21}+p_{11})\\
& &+p_{12}p_{21}p_{22}(1-p_{11}p_{12}p_{21}-p_{12}p_{21})\}[1,   2]^2\\
&=&p_{12}^{ -2}(1-p_{11}^{2}p_{12}p_{21})(1+p_{11})(p_{11}-p_{11}^{2}p_{12}p_{21}
+p_{12}^{2}p_{21}^{2})[1,   2]^2\\
&=&p_{12}^{ -2}(1-p_{11}^{2}p_{12}p_{21})(1+p_{11})(\zeta+\zeta^2\zeta+\zeta^2)[1,   2]^2=0.
\end{eqnarray*}

{\rm (e)}  $| u| =6$. $  [u] $ does not exist.
Then we show that $D=\{  [u_{4}]=[1],      [u_{1}]=[2],      [u_{2}]=[1,   2],
[u_{3}]=[1,   [1,   2]]\}$ , $p_{u_{1}
u_{2}}p_{u_{2}u_{1}}=-\zeta, p_{u_{1}u_{3}}p_{u_{3}u_{1}}=\zeta^2, p_{u_{1}u_{4}}p_{u_{4}u_{1}}=-\zeta,
p_{u_{2}u_{3}}p_{u_{3}u_{2}}=-\zeta, p_{u_{2}u_{4}}p_{u_{4}u_{2}}=-1,
 p_{u_{3}u_{4}}p_{u_{4}u_{3}}=-\zeta^2, $ and
$p_{u_{i}u_{i}}= -1,   \zeta^2,   -1,  \zeta, \ ord (p_{u_{i}u_{i}})= 2,   3,   2,  3 , \ i=1, 2, 3, 4$.
 Considering  Lemma  \ref {12}, we have
\begin{eqnarray*}
P\setminus \{1 \}
&=&\left\{u_{1},u_{2},u_{2}^2=u_2 u_2,u_{3},u_{4},u_{4}^2=u_4 u_4,u_{1}u_{2},u_{1}u_{2}^2,u_{1}u_{3},u_{1}u_{4},
u_{1}u_{4}^2,u_{2}u_{3},\right.\\
& & u_{2}^{2}u_{3},u_{2}u_{4},u_{2}^{2}u_{4}=u_2(u_2u_4),u_{2}u_{4}^2=(u_2u_4)u_4,u_{2}^{2}u_{4}^2=(u_2^2u_4)u_4,u_{3}u_{4},\\
& & u_{3}u_{4}^2,(u_{1}u_{2})u_{3},u_{1}u_{2}^{2}u_{3}=(u_{1}u_{2}^2)u_{3},(u_{1}u_{2})u_{4},(u_{1}u_{2}^{2})u_{4},
(u_{1}u_{2})u_{4}^2,\\
& &u_{1}(u_{2}^{2}u_{4}^2),u_{1}(u_{3}u_{4}),u_{1}(u_{3}u_{4}^2),(u_{2}u_{3})u_{4},(u_{2}^2u_{3})u_{4},(u_{2}u_{3})u_{4}^2,
(u_{2}^2u_{3})u_{4}^2,\\
& & \left. (u_{1}u_{2}u_{3})u_{4},u_{1}(u_{2}u_{3}u_{4}^2),u_{1}(u_{2}^2u_{3}u_{4}),u_{1}u_{2}^2u_{3}u_{4}^2
= (u_{1}u_{2})u_{2}u_{3}u_{4}^2\right\}.
\end{eqnarray*}
Thus $\mathfrak B(V) = F\oplus \mathfrak L(V)$. \end{Proof}

\begin {Proposition} \label {62} Assume that  $\mathfrak B(V) $ is connected Nichols algebra of diagonal type
and $\Delta(\mathfrak B(V)) $ is arithmetic root systems.If  $u,   v,   w \in D$
with  $\deg u  =  \deg v + \deg w$(   specially,   if $u = vw$ is the Shirshov decomposition of $u \in D$),
then $p_{wv}p_{vw}\neq1$ except the following cases:
(i) $p_{ww} = p_{vv},  p_{vv} \not =  \pm 1$;
(ii) $p_{ww} = -p_{vv}^{-1},  p_{vv} \not =  \pm 1$;
(iii) $p_{vv} = -p_{ww}^{2},  p_{ww}\in {R_{18}}$;
(iv) $p_{vv} = -p_{ww}^{-4},  p_{ww}\in {R_{18}}$;
(v) $p_{vv} = -p_{ww}^{-4},  p_{ww}\in {R_{10}}$;
(vi) $p_{ww} = -p_{vv}^{2},  p_{vv}\in {R_{18}}$;
(vii) $p_{ww} = -p_{vv}^{-4},  p_{vv}\in {R_{18}}$;
(viii) $p_{ww} = -p_{vv}^{-4},  p_{vv}\in {R_{10}}$.
\end {Proposition}

\begin{Proof} (i) $\deg (u)\in \Delta(\chi;\deg(v),  \deg(w))$ $\begin{picture}(100,        20)
\put(10,       1) {\makebox(0,  0) [t]{$\bullet$}}
\put(80,       1){\makebox(0,  0) [t]{$\bullet$}}
 \put(10,      -1) {\line(1,   0) {70}}
 \put(5,        9) {$p_{vv}$}
 \put(35,        4) {$p_{vw}p_{wv}$}
 \put(75,        9) {$p_{ww}$}
  \put(90,        1) {$.   $}\ \ \ \put(95,       1)  { }
\end{picture}$ By \cite [Prop. 2.7.1] {He05} and \cite [Lemma.2.7.2] {He05},   it is clear $p_{wv}p_{vw}\neq1$.

(ii) If exist some $k\in\mathbb N$ such that
$\deg(v)-k\deg(w)\in \mathbb N\cdot\Delta^+(\mathfrak B(V))$,   let $k_{1}\in\mathbb N$ be the maximum integer
such that $\deg(v)-k_{1}\deg(w): =  k_{2}\deg(v_{1})\in \mathbb N\cdot\Delta^+(\mathfrak B(V))$.
We know $\deg(v_{1})-k\deg(w)\notin \mathbb N\cdot\Delta(\mathfrak B(V))$ for $\ \forall k\in\mathbb N$
by the maximality of $k_{1}$. Then we obtain $\deg (u)\in \Delta(\chi;\deg(v_{1}),  \deg(w))$ with
$\begin{picture}(100,        20)
\put(10,       1) {\makebox(0,  0) [t]{$\bullet$}}
\put(80,       1){\makebox(0,  0) [t]{$\bullet$}}
 \put(10,      -1) {\line(1,   0) {70}}
 \put(5,        9) {$p_{v_{1}v_{1}}$}
 \put(35,        4) {$p_{v_{1}w}p_{wv_{1}}$}
 \put(75,        9) {$p_{ww}$}
  \put(90,        1) {$ {}   $}\ \ \ \put(95,        1)  { }
\end{picture}$ by \cite [Prop. 2.7.1] {He05} and \cite [Lemma 2.7.2] {He05}. In this case,
let $\alpha  =  \deg (v_1)$. If $\deg(v)-k\deg(w)\notin \mathbb N\cdot\Delta^+(\mathfrak B(V))$
for $\ \forall k\in\mathbb N$ and there exists some $k\in\mathbb N$ such that
$\deg(v)-k\deg(w)\in \mathbb N\cdot\Delta^-(\mathfrak B(V))$. Let $k_{1}\in\mathbb N$
be the maximum  integer such that $\deg(v)-k_{1}\deg(w): =  -k_{2}\deg(v_{1})\in \mathbb N\cdot\Delta^-(\mathfrak B(V))$.
 We know $-\deg(v_{1})-k\deg(w)\notin \mathbb N\cdot\Delta(\mathfrak B(V))$ for $\ \forall k\in\mathbb N$ by the  maximality
 of $k_{1}$.  Then we obtain $\deg (u)\in \Delta(\chi;-\deg(v_{1}),  \deg(w))$
$\begin{picture}(100,        20)
\put(10,       1) {\makebox(0,  0) [t]{$\bullet$}}
\put(80,       1){\makebox(0,  0) [t]{$\bullet$}}
 \put(10,      -1) {\line(1,   0) {70}}
 \put(5,        9) {$p_{v_{1}v_{1}}$}
 \put(35,        4) {$p_{v_{1}w}p_{wv_{1}}$}
 \put(75,        9) {$p_{ww}$}
  \put(90,        1) {${}   $}\ \ \ \put(95,        1)  { }
\end{picture}$ by \cite [Prop. 2.7.1] {He05} and \cite [Lemma 2.7.2] {He05}. In these cases,
$\deg(v) = -k_{2}\deg(v_{1})+k_{1}\deg(w)$,   $\deg(u) = -k_{2}\deg(v_{1})+(k_{1}+1)\deg(w)$,
and  $2 \nmid k_{2}$ by \cite [Cor. 2.5.4] {He05}. In this case, let $\alpha  =  -\deg (v_1)$.

\vskip.1in
 (iii) Set $\deg(w) = e_{2},  \alpha  = e_{1}$.
Then $p_{vw}p_{wv} = p_{ww}^{2k_{1}}(p_{vw}p_{wv})^{k_{2}} = p_{22}^{2k_{1}}(p_{12}p_{21})^{k_{2}}$.

T4(1).~~ $p_{0} = p_{12}p_{21}p_{11}\in R_{12}$,   $p_{11}  =  p_{0}^4$,   $p_{22}  = -p_{0}^2$,
$p_{12}p_{21} = p_{0}p_{11}^{-1} = p_{0}^{-3} = -p_{0}^{3}$.
 $p_{22}^{2k_{1}}(p_{12}p_{21})^{k_{2}} = (-p_{0}^{2})^{2k_{1}}(-p_{0}^{3})^{k_{2}} = (-1)^{k_{2}}(p_{0})^{4k_{1}+3k_{2}}\neq1$ since $2 \nmid k_{2}$.

T4(2).~~ $p_{12}p_{21}\in R_{12}$,   $p_{11}  = p_{22} = -(p_{12}p_{21})^2$,
 $p_{22}^{2k_{1}}(p_{12}p_{21})^{k_{2}} = (p_{12}p_{21})^{k_{2}}$ $(-(p_{12}p_{21})^2)^{2k_{1}} = (p_{12}p_{21})^{4k_{1}+k_{2}}\neq1$ since $2 \nmid k_{2}$.

T5(1).~~  $p_{12}p_{21}\in R_{12}$,   $p_{11}  = -(p_{12}p_{21})^2$,   $p_{22} = -1$,
 $p_{22}^{2k_{1}}(p_{12}p_{21})^{k_{2}} = (p_{12}p_{21})^{k_{2}}\neq1$ since $2 \nmid k_{2}$.

T5(2).~~  $p_{0} = p_{12}p_{21}p_{11}\in R_{12}$,   $p_{11}  = p_{0}^4$,   $p_{22} = -1$,
$p_{12}p_{21} = p_{0}p_{11}^{-1} = p_{0}^{-3} = -p_{0}^{3}$.
 $p_{22}^{2k_{1}}(p_{12}p_{21})^{k_{2}} = (-p_{0}^{3})^{k_{2}} = (-1)^{k_{2}}(p_{0})^{3k_{2}}\neq1$ since $2 \nmid k_{2}$.

T7(1).~~ $p_{11}\in R_{12}$,   $ p_{12}p_{21} = p_{11}^{-3}$,   $p_{22} = -1$;
 $p_{22}^{2k_{1}}(p_{12}p_{21})^{k_{2}} = (p_{11}^{-3})^{k_{2}} = p_{11}^{-3k_{2}}\neq1$ since $2 \nmid k_{2}$.

T7(2).~~ $p_{12}p_{21}\in R_{12}$,   $p_{11} = (p_{12}p_{21})^{-3}$,   $p_{22} = -1$.
 $p_{22}^{2k_{1}}(p_{12}p_{21})^{k_{2}} = (p_{12}p_{21})^{k_{2}}\neq1$ since $2 \nmid k_{2}$;

 T8(2)$_{1}$.~~ $ (p_{12}p_{21})^4 = -1$,   $p_{22} = -1$,   $p_{12}p_{21} = -p_{11}$;
 $p_{22}^{2k_{1}}(p_{12}p_{21})^{k_{2}} = (p_{12}p_{21})^{k_{2}}\neq1$ since $2 \nmid k_{2}$.

  T8(2)$_{2}$.~~ $ (p_{12}p_{21})^4 = -1$,   $p_{22} = -1$,   $p_{11}  = (p_{12}p_{21})^{-2}$;
 $p_{22}^{2k_{1}}(p_{12}p_{21})^{k_{2}} = (p_{12}p_{21})^{k_{2}}\neq1$ since $2 \nmid k_{2}$.

 T8(3)$.~~ (p_{12}p_{21})^4 = -1$,   $p_{11}  = (p_{12}p_{21})^2$,   $p_{22}  = (p_{12}p_{21})^{-1}$;\\
 $p_{22}^{2k_{1}}(p_{12}p_{21})^{k_{2}} = (p_{12}p_{21})^{k_{2}}(p_{12}p_{21})^{-2k_{1}} = (p_{12}p_{21})^{k_{2}-2k_{1}}\neq1$ since $2 \nmid k_{2}$.

 T10.~~  $p_{12}p_{21}\in R_{24}$,   $p_{11}  = (p_{12}p_{21})^{-6}$,   $p_{22}= (p_{12}p_{21})^{-8}$;\\
 $p_{22}^{2k_{1}}(p_{12}p_{21})^{k_{2}} = (p_{12}p_{21})^{k_{2}}((p_{12}p_{21})^{-8})^{2k_{1}} = (p_{12}p_{21})^{-16k_{1}+k_{2}}\neq1$ since $2 \nmid k_{2}$.

T11(2).~~  $p_{11}\in R_{20}$,   $p_{12}p_{21} = p_{11}^{-3}$,   $p_{22} = -1$;
 $p_{22}^{2k_{1}}(p_{12}p_{21})^{k_{2}} = (p_{11}^{-3})^{k_{2}} = p_{11}^{-3k_{2}}\neq1$ since $2 \nmid k_{2}$.

 T12.~~   $p_{11}\in{R_{30}}$,   $ p_{12}p_{21} = p_{11}^{-3}$,   $p_{22} = -p_{11}^5$;
 $p_{22}^{2k_{1}}(p_{12}p_{21})^{k_{2}} = (-p_{11}^5)^{2k_{1}}(p_{11})^{-3k_{2}}$ $= (p_{11})^{10k_{1}-3k_{2}}\neq1$ since $2 \nmid k_{2}$.

 T13.~~  $p_{12}p_{21}\in R_{24}$,   $p_{11}  = (p_{12}p_{21})^6$,   $p_{22}  = (p_{12}p_{21})^{-1}$;\\
 $p_{22}^{2k_{1}}(p_{12}p_{21})^{k_{2}} = (p_{12}p_{21})^{k_{2}}(p_{12}p_{21})^{-2k_{1}} = (p_{12}p_{21})^{-2k_{1}+k_{2}}\neq1$ since $2 \nmid k_{2}$.

 T15.~~  $p_{12}p_{21}\in R_{30}$,   $p_{11}  = -(p_{12}p_{21})^{-3}$,   $p_{22}= (p_{12}p_{21})^{-1}$;\\
 $p_{22}^{2k_{1}}(p_{12}p_{21})^{k_{2}} = (p_{12}p_{21})^{k_{2}}(p_{12}p_{21})^{-2k_{1}} = (p_{12}p_{21})^{-2k_{1}+k_{2}}\neq1$ since $2 \nmid k_{2}$.

 T16(2).~~ $p_{12}p_{21}\in R_{20}$,   $p_{11}  = (p_{12}p_{21})^{-4}$,   $p_{22} = -1$;
 $p_{22}^{2k_{1}}(p_{12}p_{21})^{k_{2}} = (p_{12}p_{21})^{k_{2}}\neq1$ since $2 \nmid k_{2}$.

 T17.~~  $p_{12}p_{21}\in R_{24}$,   $p_{11}  = -(p_{12}p_{21})^4$,   $p_{22} = -1$;
 $p_{22}^{2k_{1}}(p_{12}p_{21})^{k_{2}} = (p_{12}p_{21})^{k_{2}}\neq1$ since $2 \nmid k_{2}$.

 T18.~~  $p_{12}p_{21}\in R_{30}$,   $p_{11}  = -(p_{12}p_{21})^5$,   $p_{22} = -1$;
 $p_{22}^{2k_{1}}(p_{12}p_{21})^{k_{2}} = (p_{12}p_{21})^{k_{2}}\neq1$ since $2 \nmid k_{2}$.

 T19.~~  $p_{11}\in R_{14}$,   $ p_{12}p_{21} = p_{11}^{-3}$,   $p_{22} = -1$;
 $p_{22}^{2k_{1}}(p_{12}p_{21})^{k_{2}} = (p_{11}^{-3})^{k_{2}} = p_{11}^{-3k_{2}}\neq1$ since $2 \nmid k_{2}$.

 T20.~~  $p_{12}p_{21}\in R_{30}$,   $p_{11}  = (p_{12}p_{21})^{-6}$,   $p_{22} = -1$;
 $p_{22}^{2k_{1}}(p_{12}p_{21})^{k_{2}} = (p_{12}p_{21})^{k_{2}}\neq1$ since $2 \nmid k_{2}$.

 T21 and T22.~~  $p_{11}\in R_{24}$ or $p_{11}\in R_{14}$,   $ p_{12}p_{21} = p_{11}^{-5}$,   $p_{22} = -1$;
 $p_{22}^{2k_{1}}(p_{12}p_{21})^{k_{2}} = (p_{11}^{-5})^{k_{2}} = p_{11}^{-5k_{2}}\neq1$ since $2 \nmid k_{2}$.

T2.~~   $\Delta^+(\mathfrak B(V))  = \{e_{1},  e_{2},  e_{1}+e_{2}\}$.

T3.~~   $\Delta^+(\mathfrak B(V))  = \{e_{1},  e_{2},  e_{1}+e_{2},  2e_{1}+e_{2}\}$.

T6.~~  $p_{11}\in{R_{18}}$,   $ p_{12}p_{21} = p_{11}^{-2}$,   $p_{22} = -p_{11}^3$,
 $p_{22}^{2k_{1}}(p_{12}p_{21})^{k_{2}} = (p_{11}^{-2})^{k_{2}}(-p_{11}^{3})^{2k_{1}} = p_{11}^{6k_{1}-2k_{2}}$.\\
$\Delta^+(\mathfrak B(V))  = \{e_{1},  e_{2},  e_{1}+e_{2},  2e_{1}+e_{2},  e_{1}+2e_{2},  3e_{1}+2e_{2}\}$
by \cite [Ex. 2.5] {An}. If $\deg(v) = e_{1}+e_{2}$,   it is clear $\deg(u) = e_{1}+2e_{2},  (p_{11})^{6k_{1}-2k_{2}} = (p_{11})^{4}\neq1$,

T8(1).~~ $ p_{12}p_{21} = p_{11}^{-3}$,    $p_{22} = p_{11}^3$,   $p_{11} \in \cup _{m = 4}^\infty R_{m}$.
 $p_{22}^{2k_{1}}(p_{12}p_{21})^{k_{2}} = (p_{11}^{-3})^{k_{2}}(p_{11}^{3})^{2k_{1}} = p_{11}^{6k_{1}-3k_{2}}$.
$\Delta^+(\mathfrak B(V))  = \{e_{1},  e_{2},  e_{1}+e_{2},  2e_{1}+e_{2},  3e_{1}+e_{2},  3e_{1}+2e_{2}\}$.
If $\deg(v) = 3e_{1}+e_{2}$,   then it is clear $\deg(u) = 3e_{1}+2e_{2}$ and  $p_{11}^{6k_{1}-3k_{2}} = (p_{11})^{-3}\neq1$.

T9.~~  $p_{12}p_{21}\in R_{9}$,   $p_{11} = (p_{12}p_{21})^{-3}$,   $p_{22} = -1$;
$\Delta^+(\mathfrak B(V))  = \{e_{1},  e_{2},  e_{1}+e_{2},  2e_{1}+e_{2},  4e_{1}+3e_{2},  3e_{1}+2e_{2}\}$ by \cite [Ex. 2.5] {An},

 T11(1).~~  $p_{11}\in R_{5}$,   $ p_{12}p_{21} = p_{11}^{-3}$,   $p_{22} = -1$;
 $p_{22}^{2k_{1}}(p_{12}p_{21})^{k_{2}} = (p_{11}^{-3})^{k_{2}} = p_{11}^{-3k_{2}}$;
$\Delta^+(\mathfrak B(V))  = \{e_{1},  e_{2},  3e_{1}+e_{2},  2e_{1}+e_{2},  5e_{1}+3e_{2},  4e_{1}+3e_{2},  3e_{1}+2e_{2},  e_{1}+e_{2}\}$
by \cite [Ex. 2.7] {An}. If $\deg(v) = 3e_{1}+e_{2}$,   it is clear $\deg(u) = 3e_{1}+2e_{2}$,   then $(p_{11})^{-3k_{2}} = (p_{11})^{-9}\neq1$.

T14.~~   $p_{11}\in R_{18}$,   $ p_{12}p_{21} = p_{11}^{-4}$,   $p_{22} = -1$;
$\Delta^+(\mathfrak B(V))  = \{e_{1},  e_{2},  e_{1}+e_{2},  2e_{1}+e_{2},  4e_{1}+e_{2},  3e_{1}+e_{2}\}$ by \cite [Ex. 2.5] {An}.

T16(1).~~ $p_{11}\in R_{10}$,   $ p_{12}p_{21} = p_{11}^{-4}$,   $p_{22} = -1$;
 $p_{22}^{2k_{1}}(p_{12}p_{21})^{k_{2}} = (p_{11})^{-4k_{2}} = (p_{11})^{-4k_{2}}$.
$\Delta^+(\mathfrak B(V))  = \{e_{1},  e_{2},  3e_{1}+e_{2},  2e_{1}+e_{2},  5e_{1}+2e_{2},  4e_{1}+e_{2},  3e_{1}+2e_{2},  e_{1}+e_{2}\}$ by \cite [Ex. 2.7] {An}.
If $\deg(v) = 3e_{1}+e_{2}$,   it is clear $\deg(u) = 3e_{1}+2e_{2},  (p_{11})^{-4k_{2}} = (p_{11})^{-12}\neq1$.

\vskip.1in
(iv)  Set $\deg(w) = e_{1},  \alpha = e_{2}$. In these cases,   $\deg(v) = k_{1}e_1+k_{2}e_2$ and  $2 \nmid k_{2}$
by \cite [Cor. 2.5.4] {He05}. It is clear that $|u|\geq 3$ and $|v|\geq 2$.
Then $p_{vw}p_{wv} = p_{ww}^{2k_{1}}(p_{vw}p_{wv})^{k_{2}} = p_{11}^{2k_{1}}(p_{12}p_{21})^{k_{2}}$.

Arguments for T4 - T16 are similar to those above except for the following additional cases:

 T3(1)$_1$.~~   $p_{12}p_{21} =p_{11}^{-2}$,  $p_{22} = p_{11}^2$,  $p_{11} \not =  \pm 1$,
if $\deg(v)=e_{1}+e_{2}$, it is clear $\deg(u)=2e_{1}+e_{2}$, then $p_{11}^{2}(p_{12}p_{21})^{1}=1$.

 T3(1)$_2$.~~ $p_{12}p_{21} =p_{11}^{-2}$,  $p_{22} =-1$, $p_{11} \not =  \pm 1$,
if $\deg(v)=e_{1}+e_{2}$, it is clear $\deg(u)=2e_{1}+e_{2}$, then $p_{11}^{2}(p_{12}p_{21})^{1}=1$.

T6.~~  $p_{11}\in{R_{18}}$,   $ p_{12}p_{21}  =  p_{11}^{-2}$,   $p_{22}  =  -p_{11}^3$,
 If $\deg(v)=e_{1}+e_{2}$, it is clear $\deg(u)=2e_{1}+e_{2}$, then $(p_{11})^{2k_{1}-2k_{2}}=(p_{11})^{0}=1$.

T14.~~   $p_{11}\in R_{18}$,   $ p_{12}p_{21}  =  p_{11}^{-4}$,   $p_{22}  =  -1$;
 If $\deg(v)=2e_{1}+e_{2}$, it is clear $\deg(u)=3e_{1}+e_{2}$, then $(p_{11})^{2k_{1}-4k_{2}}=(p_{11})^{0}=1$.

T16(1).~~ $p_{11}\in R_{10}$,   $ p_{12}p_{21}  =  p_{11}^{-4}$,   $p_{22}  =  -1$;
If $\deg(v)=2e_{1}+e_{2}$, it is clear $\deg(u)=3e_{1}+e_{2}$, then $(p_{11})^{2k_{1}-4k_{2}}=(p_{11})^{0}=1$.

\vskip.1in
Similarly, we have (v)--(viii).  \end{Proof}

\begin{Remark}
(i) If $[u][v]\neq 0$, then the following conditions are equivalent:
(1) $[[u], [v]]=0$ ,  $[[v], [u]]=0$;
(2) $1-p_{uv}p_{vu}=0$ ,  $[[v], [u]]=0$;
(3) $1-p_{uv}p_{vu}=0$ ,  $[[u], [v]]=0$.

(ii) If $[[u], [v]]=0$,  then $[[v], [u]]=(1-p_{uv}p_{vu})[u][v]$.
\end{Remark}

\begin{Remark} There is a braided m-Lie algebra which is not a  Nichols braided Lie algebra
or Nichols braided m-Lie algebra. In fact,   let $L = L_{\bar 0} + L_{\bar 1}$ be a super-Lie
algebra with    $L_{\bar 0} = sl (2)$,   $ L_{\bar 1} =0$.
It is clear that $L$ is a finite dimensional  m-braided Lie algebra of  diagonal type
in $^{k \mathbb Z_2}_{k \mathbb Z_2} { \mathcal  YD}$.
Because $L$ is not nilpotent and every finite dimensional $\mathfrak L(V)$ is nilpotent,
$L$ is not a Nichols braided Lie algebra or Nichols braided m-Lie algebra.
\end{Remark}

\section {Cartan type} \label {s5}

In this section  we give the sufficient conditions for Nichols braided Lie algebra $\mathfrak L(V)$ 
to be a homomorphic image of  a braided Lie algebra generated by  $V$ with defining relations. 
Basic field $F$ is the complex field $\mathbb C.$

Let $\Phi^+$ denote the positive root system  of simple Lie algebras. $E_e = n-1, n-1, n-1, n-1, 5, 6, 7, 3, 1$
in $A_n, B_n, C_n, D_n, E_6, E_7, E_8, F_4, G_2$, respectively, i.e. $E_e = {\rm rank } \Phi -1$. In fact,  $\Phi^+ = \Delta ^+(\mathfrak B(V))$
by \cite {He06b}. Let $\epsilon _1,     \epsilon _2,     \cdots,     \epsilon_n$ be a normal orthogonal basis of
$\mathbb R^n$;  $I   =  : \sum _{i  =  1}^n\mathbb Z
\epsilon _i$ and $I'   =  I + \mathbb Z (\epsilon _1+ \epsilon _2 +\cdots + \epsilon_n)/2.   $
Let $\Delta   =   \{e_1, \cdots, e_n\}$ be a prime root system.

\begin {Lemma} \label {4.1}
Let $\Omega$ be the class of connected components  of $V$. then
$\mathfrak{B} (V) = \oplus _{J\in \Omega }\mathfrak{B}(V_J)$  and
$\mathfrak{L} (V) = \oplus _{J\in \Omega }\mathfrak{L}(V_J).$
\end {Lemma}

\begin {Lemma} \label {4.1''}
Assume that $V$ is a braided vector space of   diagonal type with braided matrix $(q_{ij})_{n\times n}$
and basis $x_1, x_2,  \cdots,  x_n$. Then ${\rm ord } ( q_{ij})$ is finite for any $1\le i,  j \le n$
if and only if there exists a finite abelian group $G$ such that  $V$ becomes a $kG$- {\rm YD} module.
\end {Lemma}

\begin{Proof} The necessity  is clear. The sufficiency. Let $N$ be the  least common multiple
  of $\{ ord q_{ij} \mid 1\le i,  j \le n \}$ and $G = (g_1) \times (g_2) \times \cdots \times (g_n)$
with ${\rm ord } ( g_i) =N$ for $1\le i\le n$. Let $\chi (g_1 ^{k_1} \cdots g_n ^{k_n},
g_1 ^{k_1'} \cdots g_n ^{k_n'} ):= \prod _{1\le i,  j \le n}q_{ij} ^{k_ik_j'}$.
It is clear that $\chi $ is a bicharacter on $G \times G$.
Therefore $V$ becomes a $kG$- {\rm YD} module. \end{Proof}

\begin {Lemma} \label {4.1'} Assume that $V$ is a braided vector space of   diagonal type with
braided matrix $(q_{ij})_{n\times n}$.
If  $\dim  \mathfrak B(V)<\infty$, then   there exists a braided matrix $(q_{ij}')_{n \times n}$,
which is twisted equivalent to $(q_{ij})_{n\times n}$,  such that  ${\rm ord } ( q_{ij}' ) <\infty$
for any $1\le i,  j \le n.$
\end {Lemma}

\begin{Proof} We show this by two steps.

(i) ${\rm ord } ( q_{ij}q_{ji}) < \infty$ for any $1\le i, j \le n$.
In fact, it is clear ${\rm ord } (q_{ii} ^2 )< \infty$ for any $1\le i \le n.$
If there exist $i$ and $j$ with $i < j$ such that ${\rm ord } ( q_{ij}q_{ji}) = \infty$.
Obviously, $[u]:= [x_i, x_j] \in D.$ Consequently, ${\rm ord } ( p_{u, u}) =\infty$ and
$\dim  \mathfrak B(V)=\infty$, which is a contradiction.
(ii) Set $q_{ij}' := \sqrt{q_{ij} q_{ji}}$ for any $1\le i, j \le n$. \end{Proof}

\begin {Lemma} \label {4.222'} Assume that $(V,   (q_{ij})_{n\times n})$ is of connected
Cartan type  with Cartan matrix $(a_{ij})_{n\times n}$. Then

(i) For $A_n, D_n, E_8, E_7, E_6$:
 $p_{\alpha, \alpha}   =     q $, $p_{\alpha, \gamma}p_{ \gamma, \alpha} = q^{(\alpha, \gamma)}$,
$p_{\alpha, \beta} p_{ \beta, \alpha} = q^{\pm 1}$ or $1$
for any $\alpha, \beta,   \gamma\in \Delta^+(\mathfrak B(V))$ with $\alpha \not= \beta$.

(ii) For $B_n$:

\ \ \ \ \ \ \ \ \  \ \ \
$\begin{picture}(100,     15)
\put(27,     1){\makebox(0,    0)[t]{$\bullet$}}
\put(60,     1){\makebox(0,     0)[t]{$\bullet$}}
\put(93,    1){\makebox(0,    0)[t]{$\bullet$}}
\put(159,     1){\makebox(0,     0)[t]{$\bullet$}}
\put(192,    1){\makebox(0,     0)[t]{$\bullet$}}
\put(225,    1){\makebox(0,    0)[t]{$\bullet$}}
\put(28,     -1){\line(1,     0){33}}
\put(61,     -1){\line(1,     0){30}}
\put(130,    -1){\makebox(0,    0)[t]{$\cdots\cdots\cdots\cdots$}}
\put(160,    -1){\line(1,     0){30}}
\put(193,     -1){\line(1,     0){30}}
\put(22,    -15){1}
\put(58,     -15){2}
\put(91,     -15){3}
\put(157,     -15){n-2}
\put(191,     -15){n-1}
\put(224,     -15){n}
\put(22,    10){$q^{2}$}
\put(58,     10){$q^{2}$}
\put(91,     10){$q^{2}$}
\put(157,     10){$q^{2}$}
\put(191,     10){$q^{2}$}
\put(224,     10){$q$}
\put(40,     5){$q^{-2}$}
\put(73,     5){$q^{-2}$}
\put(172,    5){$q^{-2}$}
\put(205,     5){$q^{-2}$}
\end{picture}$
\vskip.1in
\begin{eqnarray*}
\Delta^+(\mathfrak B(V))  &=& \left \{\sum \limits_{i  =  n}^{n}e_{i},  \sum \limits_{i  =  n-1}^{n}e_{i},
  \sum \limits_{i  =  n-2}^{n}e_{i},  \cdots,  \sum \limits_{i  =  1}^{n}e_{i}\right \}\\
& & \cup ~\left \{ \sum \limits_{i  =  1}^{1}e_{i}+\sum \limits_{i  =  2}^{n}2e_{i},
\sum \limits_{i  =  1}^{2}e_{i}+\sum \limits_{i  =  3}^{n}2e_{i},  \sum \limits_{i  =  1}^{3}e_{i}+\sum \limits_{i  =  4}^{n}2e_{i}, \right. \\
& &  \cdots,\sum \limits_{i  =  1}^{n-1}e_{i}+\sum \limits_{i  =  n}^{n}2e_{i},
\sum \limits_{i  =  2}^{2}e_{i}+\sum \limits_{i  =  3}^{n}2e_{i},  \sum \limits_{i  =  2}^{3}e_{i}+\sum \limits_{i  =  4}^{n}2e_{i},    \\
& & \left.\cdots,\sum \limits_{i  =  2}^{n-1}e_{i}+\sum \limits_{i  =  n}^{n}2e_{i},  \cdots,
\sum \limits_{i  =  n-1}^{n-1}e_{i}+\sum \limits_{i  =  n}^{n}2e_{i} \right\}\\
& &\cup~ \left\{ \sum \limits_{i  =  1}^{1}e_{i},  \sum \limits_{i  =  1}^{2}e_{i},
\sum \limits_{i  =  1}^{3}e_{i},  \cdots,  \sum \limits_{i  =  1}^{n-1}e_{i},  \sum \limits_{i  =  2}^{2}e_{i},  \sum \limits_{i  =  2}^{3}e_{i},    \right.\\
& & \left. \cdots,\sum \limits_{i  =  2}^{n-1}e_{i},  \cdots, \sum \limits_{i  =  n-2}^{n-2}e_{i},
\sum \limits_{i  =  n-2}^{n-1}e_{i},  \sum \limits_{i  =  n-1}^{n-1}e_{i}\right\}:  =  Q \cup  S\cup  T,
\end{eqnarray*}
$p_{\alpha,   \alpha }  =  q$ for $\ \forall \alpha \in Q$ and
$p_{\alpha,   \alpha}  =  q^2$ for $\ \forall \alpha \in S$ or $T$.

(iii) For  $C_n$:

\ \ \ \ \ \ \ \ \  \ \ \
$\begin{picture}(100,     15)
\put(27,     1){\makebox(0,    0)[t]{$\bullet$}}
\put(60,     1){\makebox(0,     0)[t]{$\bullet$}}
\put(93,    1){\makebox(0,    0)[t]{$\bullet$}}
\put(159,     1){\makebox(0,     0)[t]{$\bullet$}}
\put(192,    1){\makebox(0,     0)[t]{$\bullet$}}
\put(225,    1){\makebox(0,    0)[t]{$\bullet$}}
\put(28,     -1){\line(1,     0){33}}
\put(61,     -1){\line(1,     0){30}}
\put(130,    -1){\makebox(0,    0)[t]{$\cdots\cdots\cdots\cdots$}}
\put(160,    -1){\line(1,     0){30}}
\put(193,     -1){\line(1,     0){30}}
\put(22,    -15){1}
\put(58,     -15){2}
\put(91,     -15){3}
\put(157,     -15){n-2}
\put(191,     -15){n-1}
\put(224,     -15){n}
\put(22,    10){$q$}
\put(58,     10){$q$}
\put(91,     10){$q$}
\put(157,     10){$q$}
\put(191,     10){$q$}
\put(224,     10){$q^2$}
\put(40,     5){$q^{-1}$}
\put(73,     5){$q^{-1}$}
\put(172,    5){$q^{-1}$}
\put(205,     5){$q^{-2}$}
\end{picture}$
\vskip.1in
\begin{eqnarray*}
\Delta^+(\mathfrak B(V))   &=& \left \{\sum \limits_{i  =  1}^{n-1}2e_{i}+e_{n},
\sum \limits_{i  =  2}^{n-1}2e_{i}+e_{n},  \cdots,  \sum \limits_{i  =  n-1}^{n-1}2e_{i}+e_{n},  e_{n}\right \}\\
& & \cup ~ \left\{\sum \limits_{i  =  1}^{1}e_{i}+\sum \limits_{i =  2}^{n-1}2e_{i}+e_{n},
\sum \limits_{i  =  1}^{2}e_{i}+\sum \limits_{i  =  3}^{n-1}2e_{i}+e_{n}, \right.\\
& &   \sum \limits_{i =  1}^{3}e_{i}+\sum \limits_{i  =  4}^{n-1}2e_{i}+e_{n}, \cdots,
\sum \limits_{i  =  1}^{n-1}e_{i}+e_{n},  \sum \limits_{i  =  2}^{2}e_{i}+\sum \limits_{i  =  3}^{n-1}2e_{i}+e_{n},   \\
& & \left. \sum \limits_{i  =  2}^{3}e_{i}+
\sum \limits_{i  =  4}^{n-1}2e_{i}+e_{n}, \cdots, \sum \limits_{i  =  2}^{n-1}e_{i}+e_{n},  \cdots,
\sum \limits_{i  =  n-1}^{n-1}e_{i}+e_{n}\right \} \\
& & \cup ~ \left\{\sum \limits_{i  =  1}^{1}e_{i},  \sum \limits_{i  =  1}^{2}e_{i},
\sum \limits_{i  =  1}^{3}e_{i},  \cdots,  \sum \limits_{i  =  1}^{n-1}e_{i},  \sum \limits_{i  =  2}^{2}e_{i},  \sum \limits_{i  =  2}^{3}e_{i}, \right. \\
& &\left.  \cdots,\sum \limits_{i  =  2}^{n-1}e_{i},  \cdots,  \sum \limits_{i  =  n-2}^{n-2}e_{i},
\sum \limits_{i  =  n-2}^{n-1}e_{i},  \sum \limits_{i  =  n-1}^{n-1}e_{i}\right \} :  =  Q\cup  S\cup T,
\end{eqnarray*}
$p_{\alpha,   \alpha }  =  q^2$ for $\ \forall \alpha \in Q$ and
$p_{\alpha,   \alpha}  =  q$ for $\ \forall \alpha \in S$ or $T$.

(iv) For $F_4$:

 \ \ \ \ \ \ \ \ \  \ \ \
$\begin{picture}(100,     15)
\put(27,     1){\makebox(0,    0)[t]{$\bullet$}}
\put(60,     1){\makebox(0,     0)[t]{$\bullet$}}
\put(93,    1){\makebox(0,    0)[t]{$\bullet$}}
\put(126,     1){\makebox(0,   0)[t]{$\bullet$}}
\put(28,     -1){\line(1,     0){33}}
\put(61,     -1){\line(1,     0){30}}
\put(94,    -1){\line(1,     0){30}}
\put(22,    -15){1}
\put(58,     -15){2}
\put(91,     -15){3}
\put(124,    -15){4}
\put(22,    10){$q^2$}
\put(58,     10){$q^2$}
\put(91,     10){$q$}
\put(124,     10){$q$}
\put(40,     5){$q^{-2}$}
\put(73,     5){$q^{-2}$}
\put(106,    5){$q^{-1}$}
\end{picture}$
\vskip.1in
\begin{eqnarray*}
\Delta^+(\mathfrak B(V))  & = &  \{ e_2+2 e_3+ 2 e_4,   e_1+ e_2+2 e_3+ 2 e_4, e_1+ 2e_2+2 e_3+ 2 e_4,\\
& & e_1,e_1+e_2,e_2, 2e_1+ 3e_2+4 e_3+ 2 e_4,e_1+ 3e_2+4 e_3+ 2 e_4,\\
& & e_1+ 2e_2+4 e_3+ 2 e_4,e_1+ 2e_2+2e_3,e_1 + e_2+2e_3,e_2+2e_3\}\\
& & \cup ~\{e_1+ 2e_2+3 e_3+  e_4,   e_2+2 e_3+  e_4,e_1+ e_2+2 e_3+  e_4,\\
& &e_1+ 2e_2+2e_3+ e_4,e_3+  e_4,e_2+ e_3+  e_4,e_1+ e_2+ e_3+ e_4,e_4,\\
& & e_1+ 2e_2+3 e_3+ 2 e_4,e_1+ e_2+ e_3,e_2+  e_3,e_3  \} :  =   Q\cup S,
\end{eqnarray*}
 $p_{\alpha,    \alpha }  =  q^2$ for $\ \forall\alpha \in Q $ and  $p_{\alpha,    \alpha }
=  q$ for $\ \forall\alpha \in S $.

(v) For $G_2$:
\vskip.1in
\ \ \ \ \ \ \ \ \  \ \ \
$\begin{picture}(100,     15)
\put(27,     1){\makebox(0,    0)[t]{$\bullet$}}
\put(60,     1){\makebox(0,     0)[t]{$\bullet$}}
\put(28,     -1){\line(1,     0){33}}
\put(22,    -15){1}
\put(58,     -15){2}
\put(22,    10){$q$}
\put(58,     10){$q^3$}
\put(40,     5){$q^{-3}$}
\end{picture}$\\ \\
$\Delta^+(\mathfrak B(V))   =  \{e_{1},  e_{1}+e_{2},  2e_{1}+e_{2}\} \cup \{ 3e_{1}+e_{2},  3e_{1}+2e_{2},  e_{2}\}:  =  Q \cup S.$
$p_{\alpha\alpha}  =  q$ for $\ \forall \alpha\in Q$ and
$p_{\alpha\alpha}  =  q^3$ for $\ \forall \alpha \in S$.
\end {Lemma}

\begin{Proof} (i)  By \cite [Section 12.1]{Hu72}, the root system
$\Delta^+(\mathfrak B(V))  =   \{ \beta \in I \mid (\beta, \beta )   =  2 \}.$
Let $\alpha   =   k_1e_1 + \cdots + k _n e_n, \gamma   =   k_1'e_1 + \cdots + k _n' e_n    \in \Delta^+(\mathfrak B(V))$.
\begin {eqnarray}  \label {e6.1.2}   (\alpha, \gamma )  =  \sum _{i  =  1}^n 2k_ik'_i + \sum _{i \not= j} a_{ij}k_ik_j'
\end {eqnarray} and
\begin {eqnarray}  \label {e6.1.3}
p_{\alpha, \gamma }p _{\gamma, \alpha}  = q ^ { \sum _{i  =  1}^n 2k_ik'_i + \sum _{i \not= j} a_{ij}k_ik_j'} = q ^{(\alpha, \gamma)}.
\end {eqnarray}
 Consequently,
 $p_{\alpha, \alpha}   =    q $ by (\ref {e6.1.3}). By \cite [Section 9.4, Table 1]{Hu72},  $(\alpha, \beta ) = 1 $ or $-1$ or $0.$

 (ii) - (v) are clear. \end{Proof}

\vskip.1in
Recall that $(T(V), [\ ])$ is a braided Lie algebra. The braided Lie algebra  $(FL(V), [\ ])$
generated by $V$ in $(T(V), [\ ])$ is called the free braided Lie algebra of $V.$
If $f_1, f_2, \cdots, f_r \in FL(V)$ and $I$ is an ideal $I$ generated by $f_1, f_2, \cdots, f_r$
in $(FL(V), [\ ])$, then $(FL(V)/ I, [\ ])$ is called the braided Lie algebra generated by
$x_1, x_2, \cdots, x_n$ with the defining relations  $f_1, f_2, \cdots, f_r$.

\begin {Lemma} \label {4.2'} Assume that $(V, (q_{ij})_{n\times n})$ is a connected braided vector space
of finite Cartan type  with Cartan matrix $(a_{ij})_{n\times n}$. If  ${\rm ord } (q_{ii})$
is prime to $3$ when $q_{ij} q_{ji} \in \{ q_{ii}^ {3},  q_{jj} ^{3}\}$,
then  ${\rm ord } ( p_{\alpha, \alpha } )=N$  for  root vector $x_\alpha$
with $\alpha \in \Delta^+(\mathfrak B(V))$ and $N = ord (q_{11})$,
where root vectors were defined in \cite {Lu90}.
\end {Lemma}

\begin{Proof} By \cite [Th.1.1(i)]{AS00}, ${\rm ord } ( q_{ii}) =N$ for $1\le i \le n.$
 $ord (p_{\alpha, \alpha })=  N$ for any root $\alpha$ by Lemma \ref {4.222'} and
  $x_\alpha^i \in \mathfrak L(V) $ for $1\le i \le N$  by Lemma \ref {11}. \end{Proof}

 \vskip.1in
Let $ad _c xy := [x, y]_c = xy - p_{x, y} yx.$
\begin {Theorem} \label {4.3} If
$V$ is a finite Cartan type  with Cartan matrix $(a_{ij})_{n\times n}$ and  the following
conditions satisfied for any $1\le i, j \le n:$
(i) ${\rm ord }  (q_{ii})$ is odd;
(ii)  ${\rm ord }  (q_{ii})$ is prime to $3$ when $q_{ij} q_{ji} \in \{ q_{ii}^ {3},  q_{jj} ^{3}\};$
(iii) ${\rm ord } ( q_{ij})< \infty$.
then Nichols braided Lie algebra $\mathfrak{L} (V)$ is a homomorphic image of the braided Lie algebra 
generated by $x_1,    x_2,    \cdots,    x_n$ with defining relations:
(iv)  $ad _c x_i ^{1- a_{ij}}x_j$,    $i\not  =   j.$
(v) $x_\alpha ^N$ for any $\alpha \in \Delta ^+(\mathfrak B(V)),   $ where $N$ is order of $q_{11}$.
\end {Theorem}

\begin{Proof} By \cite [Section 4.1] {AS02},     $x_\alpha \in FL(V)$. It follows from Lemma \ref {4.2'} and Lemma
\ref {11} that $x_\alpha ^N \in FL(V)$. Let $I$ and $J$ denote ideals generated by elements of  (iv) and (v) in
$T(V)$ as algebras and in $FL(V)$ as braided Lie algebras with bracket  $[\ \ ]$. Consequently,   using  \cite [Theorem 5.1] {AS10},
we have that the map from $FL(V)/ J$ to $\mathfrak L(V)$  by sending $x +J$ to $x+I$ is a epimorphism. \end{Proof}

\begin {Theorem} \label {4.2}  If $\dim \mathfrak{B} (V)< \infty$ and the following conditions satisfied for any $1\le i, j \le n:$
(i) ${\rm ord }  (q_{ii})$ is odd;
(ii)  ${\rm ord }  (q_{ii})$ is prime to $3$ when $q_{ij} q_{ji} \in \{ q_{ii}^ {3},  q_{jj} ^{3}\};$
(iii) ${\rm ord } (q_{ii})>3$;
(iv)  ${\rm ord } ( q_{ij})< \infty$.
Then $V$ is a finite Cartan type  with Cartan matrix $(a_{ij})_{n\times n}$ and
Nichols braided Lie algebra $\mathfrak{L} (V)$ is  a homomorphic image of the braided Lie algebra generated  by $x_1,  x_2,  \cdots,  x_n$
with the defining relation,
(v)  $ad _c x_i ^{1- a_{ij}}x_j$,  $i\not= j$;
(vi)  $x_\alpha ^N$ for any $\alpha \in \Delta ^+(\mathfrak B(V)), $ where $N$ is order of $q_{11}$.
\end {Theorem}

\begin{Proof}  By the proof of  \cite [Theorem 5.3] {AS10}, $V$ is a finite Cartan type.
Using Theorem \ref {4.3} we complete the proof. \end{Proof}

\begin {Lemma} \label {8''}  Under the conditions of Theorem \ref {4.2} or Theorem \ref {4.3},
for any $\alpha \in\Delta ^+(\mathfrak B(V))$, there exists a unique hard super-letter $[u]$
such that $\deg (u)   =   \alpha$.
 \end {Lemma}

\begin{Proof}  Considering the dimensional formulas of $\mathfrak B(V)$ in
Theorem \ref {8} and \cite [Th. 1.1(i)]{AS00}, we complete the proof. \end{Proof}

\begin {Lemma} \label {8'''}   Assume  that  $(V,   (q_{ij})_{n\times n})$ is  of  a  connected
Cartan type.If $\alpha, \beta, \gamma \in\Delta ^+(\mathfrak B(V))$ with $\alpha + \beta = \gamma$,
 then $p_{\alpha, \beta} p_{ \beta, \alpha} =1 $ if and only if  $(\alpha, \beta)$
or $( \beta, \alpha) \in X$  and $X$ is defined in the following cases:

(i) For $F_4$, $X := \{  (\epsilon _i, \epsilon _j) \mid $ $i\not= j \}$
$ \cup ~ \{ (\frac {1}{2}  (\epsilon_1 + (-1)^{k_2} \epsilon_2+ (-1)^{k_3} \epsilon_3+ (-1)^{k_4}\epsilon_4 ),
 \frac {1}{2}  (\epsilon_1 + (-1)^{k_2'} \epsilon_2+ (-1)^{k_3'} \epsilon_3+ (-1)^{k_4'}\epsilon_4 ) )
 \mid   (-1)^{k_2}+ (-1)^{k_2'} ~{\rm or} ~ (-1)^{k_3}+(-1)^{k_3'}~{\rm or}~
 (-1)^{k_4}+(-1)^{k_4'}~  {\rm is~ } 1  \}$;

(ii)  For $B_n$, $X= \{ (\epsilon _i, \epsilon _j \mid 1\le i \not= j \le n\}$;

(iii)  For $C_n$,  $X:= \{ (\epsilon _i -\epsilon _j, \epsilon _i+\epsilon _j)\mid  1\le i< j \le n \}$.
 \end {Lemma}

\begin{Proof} By Proposition \ref {62}, it is clear that   $p_{\alpha, \beta} p_{ \beta, \alpha} =1$  if and only if
 \begin {eqnarray} \label {e5.1}
p_{\gamma, \gamma} = q^2, p_{\alpha, \alpha} = p_{\beta, \beta}=q;
\end {eqnarray}
 \begin {eqnarray} \label {e5.2}
\hbox { or  } \ \ \  p_{\gamma, \gamma} = q,  p_{\alpha, \alpha} = p_{\beta, \beta}=q^2, q^3 =1.
\end {eqnarray}
By Lemma \ref {4.2'}, $ p_{\alpha, \alpha} =q$ and  $p_{\alpha, \beta} p_{ \beta, \alpha} \not=1 $
for  $A_n, D_n, E_8, E_7, E_6$. We can complete the proof by simple computation. \end{Proof}

\begin {Proposition} \label {5.6} Assume  $[u]\in D$.
(i)  If  $p_{vw}p_{wv} \not= 1$ for any two  descendants $v$ and $w$ of $u$, then $[u]^- \in {\mathfrak L}(V)$;
(ii) If  connected $(V,   (q_{ij})_{n\times n})$ is  of $A_n, D_n, E_6, E_7, E_8, G_2$, then $[u]^- \in {\mathfrak L}(V)$.
\end {Proposition}

\begin{Proof}  (i) By induction  and Lemma \ref {12}, we obtain (i).
(ii) follows from (i) and Lemma \ref {8'''}.  \end{Proof}

\begin {Theorem} \label {5.6'} Under the conditions of Theorem \ref {4.2} or Theorem \ref {4.3},
assume that   connected $(V,   (q_{ij})_{n\times n})$ is  of $A_n, D_n, E_6, E_7, E_8, G_2$,
and $q= q_{ii}$ with $N := {\rm ord } (q_{11})$.
(i) If  $u, v\in D,$  then $p_{uu}^{i} p_{u, uv}p_{uv, u}\not= 1$  for $1\le i \le  2([\frac {N}{2}]-1)-2$.
(ii) $ \dim {\mathfrak L}(V) \ge  (N-1)^{ \mid \Phi^+ \mid } + ( [\frac {N} {2}]-1)\frac { \mid \Phi^+ \mid  ( \mid \Phi^+ \mid-1  )}{2}$.
\end {Theorem}

\begin{Proof}   (i) By Lemma \ref {4.2'} (i),
$p_{uu}^{i} p_{u, uv}p_{uv, u} = p_{uu} ^ip_{uu}^2 p_{uv}p_{vu}= q ^{i +2+ ( \deg (u), \deg (v))} \not= 1 $
for $1\le i \le  2([\frac {N}{2}]-1)-2$.
(ii)  Let $D:= \{ u_1, u_2, \cdots, u_r\}$ with $u_r < u_{r-1} < \cdots < u_1 $.
By Lemma \ref {8''}, $r = \mid \Phi ^+ \mid.$ It follows from Lemma \ref {21} that $u^i uv \in \mathfrak L(V)$
for any $u, v\in D$ and $1\le i \le  [\frac{N}{2}] -1. $ Consequently,
$B:= \{  u_1^j, u_2^j, $ $\cdots, u_r^j;  u_1^iu_1u_2, u_1^iu_1u_3, $ $ \cdots,
u_1^iu_1u_r;  u_2^iu_2u_3, u_2^iu_2u_4, \cdots,  u_2^iu_2u_r; $ $ \cdots ;  $ $  u_{r-1}^iu_{r-1}u_r \} $ $\subset P \cap {\mathfrak L}(V) $,
which implies
 $\dim {\mathfrak L}(V) \ge  (N-1)^{ \mid \Phi^+ \mid } + ( [\frac {N} {2}]-1)\frac { \mid \Phi^+ \mid  ( \mid \Phi^+ \mid-1  )}{2}$. \end{Proof}

\vskip.3in
\noindent {\bf\large Acknowledgments}
\vskip.1in
\noindent YZZ acknowledges the partial financial support from the Australian Research Council through
Discovery Projects DP110103434 and DP140101492. We would like to thank the referee for many suggestions
which lead to the substantial improvement of the paper.

\end {document}